\documentclass[12pt]{article}

\usepackage{amsmath, amssymb, amscd, fancyhdr}

\newcommand{\version}{version 1.0,\ \ Nov. 25, 2012}

\def\eqref#1{(\ref{#1})}
\newcommand{\goth}{\mathfrak}
\newcommand{\g}{{\frak g}}
\newcommand{\arrow}{{\:\longrightarrow\:}}

\newcommand{\C}{{\Bbb C}}
\newcommand{\R}{{\Bbb R}}

\newcommand{\6}{\partial}
\def\1{\sqrt{-1}\:}
\newcommand{\restrict}[1]{{\left|_{{\phantom{|}\!\!}_{#1}}\right.}}
\newcommand{\cntrct}                
{\hspace{2pt}\raisebox{1pt}{\text{$\lrcorner$}}\hspace{2pt}}
\newcommand{\3}{{/\!\!/\!\!/}} 
\newcommand{\2}{{/\!\!/}}

\def\Bbb#1{\mathbb #1}

\newcommand{\calo}{{\cal O}}

\renewcommand{\phi}{\varphi}
\renewcommand{\epsilon}{\varepsilon}
\renewcommand{\geq}{\geqslant}
\renewcommand{\leq}{\leqslant}

\newcommand{\End}{\operatorname{End}}
\newcommand{\Sym}{\operatorname{Sym}}
\newcommand{\Sp}{\operatorname{Sp}}
\newcommand{\Mer}{\operatorname{Mer}}
\newcommand{\Tot}{\operatorname{Tot}}
\newcommand{\Id}{\operatorname{Id}}

\newcommand{\Hom}{\operatorname{Hom}}

\newcommand{\Lie}{\operatorname{Lie}}
\newcommand{\Tw}{\operatorname{Tw}}

\newcommand{\comment}[1]{{}}

\def\endproof{\blacksquare}

\makeatletter

\newcounter{Mycounter}[section]
\newcounter{lemma}[section]
\renewcommand{\thelemma}{\noindent{Lemma \thesection.\arabic{lemma}}}
\newcommand{\lemma}{%
     \setcounter{lemma}{\value{Mycounter}}
     \refstepcounter{lemma}
     \stepcounter{Mycounter}
     {\bf \thelemma:\ }}

\newcounter{claim}[section]
\renewcommand{\theclaim}{\noindent{Claim \thesection.\arabic{claim}}}
\newcommand{\claim}{%
     \setcounter{claim}{\value{Mycounter}}
     \refstepcounter{claim}
     \stepcounter{Mycounter}
     {\bf \theclaim:\ }}

\newcounter{sublemma}[section]
\newcounter{corollary}[section]
\renewcommand{\thecorollary}{\noindent{Corollary \thesection.\arabic{corollary}}}
\newcommand{\corollary}{%
     \setcounter{corollary}{\value{Mycounter}}
     \refstepcounter{corollary}
     \stepcounter{Mycounter}
     {\bf \thecorollary:\ }}

\newcounter{theorem}[section]
\renewcommand{\thetheorem}{\noindent{Theorem \thesection.\arabic{theorem}}}
\newcommand{\theorem}{%
     \setcounter{theorem}{\value{Mycounter}}
     \refstepcounter{theorem}
     \stepcounter{Mycounter}
     {\bf \thetheorem:\ }}

\newcounter{conjecture}[section]
\newcounter{proposition}[section]
\renewcommand{\theproposition}
       {\noindent{Proposition \thesection.\arabic{proposition}}}
\newcommand{\proposition}{%
     \setcounter{proposition}{\value{Mycounter}}
     \refstepcounter{proposition}
     \stepcounter{Mycounter}
     {\bf \theproposition:\ }}

\newcounter{definition}[section]
\renewcommand{\thedefinition}
       {\noindent{Definition~\thesection.\arabic{definition}}}
\newcommand{\definition}{%
     \setcounter{definition}{\value{Mycounter}}
     \refstepcounter{definition}
     \stepcounter{Mycounter}
     {\bf \thedefinition:\ }}

\newcounter{example}[section]
\newcounter{remark}[section]
\renewcommand{\theremark}{\noindent{Remark \thesection.\arabic{remark}}}
\newcommand{\remark}{%
     \setcounter{remark}{\value{Mycounter}}
     \refstepcounter{remark}
     \stepcounter{Mycounter}
     {\bf \theremark:\ }}

\newcounter{problem}[section]
\newcounter{question}[section]
\@addtoreset{equation}{section}
\@addtoreset{footnote}{section}
\makeatother

\begin{document}

\begin{center}
{\LARGE\bf Holography principle for twistor spaces
}
\\[4mm]
Misha Verbitsky\footnote{Partially supported by RFBR grants
 12-01-00944-Á,  10-01-93113-NCNIL-a, and
AG Laboratory NRI-HSE, RF government grant, ag. 11.G34.31.0023.}
\\[4mm]

{\tt verbit@verbit.ru}

\end{center}

{\small 
\hspace{0.15\linewidth}
\begin{minipage}[t]{0.7\linewidth}
{\bf Abstract} \\
Let $S$ be a smooth rational curve on a complex manifold $M$.
It is called {\bf ample} if its normal bundle is positive:
$NS=\bigoplus\calo(i_k)$, $i_k>0$. We assume that $M$ is
covered by smooth holomorphic deformations of $S$.
The basic example of such a manifold is a twistor space
of a hyperk\"ahler or a 4-dimensional anti-selfdual
Riemannian manifold $X$ (not necessarily compact).
We prove ``a holography principle'' for such a manifold:
any meromorphic function defined in a neighbourhood $U$
of $S$ can be extended to $M$, and any section of a 
holomorphic line bundle can be extended from $U$ to $M$.
This is used to define the notion of a {\bf Moishezon
  twistor space}: this is a twistor space admitting a holomorphic
embedding to a Moishezon variety $M'$. We show that 
this property is local on $X$, and the 
variety $M'$ is unique up to birational transform. 
We prove that the twistor spaces of hyperk\"ahler manifolds
obtained by hyperk\"ahler reduction of flat 
quaternionic-Hermitian spaces by the action 
of reductive Lie groups (such as Nakajima's
quiver varieties) are always Moishezon.
\end{minipage}
}

{
\small
\tableofcontents
}


\section{Introduction}


\subsection{Quasilines on complex manifolds}

The present paper was written as an attempt to answer the following
question. Let $S\subset M$ be a smooth rational curve
in a complex manifold, with normal bundle $NS$ isomorphic
to $\calo(1)^n$.\footnote{Such rational curves are called {\bf quasilines}, 
see e.g. \cite{_Badescu_Beltrametti_Ionescu_}.}
Is there a notion of a normal form for a tubular
neighbourhood of such a curve?

When the normal bundle is $\calo(-1)^n$ instead of $\calo(1)^n$, 
a tubular neighbourhood  of the curve has a normal form,
obtained by blowing down this curve to a point, and taking
a sufficiently small Stein neighbourhood of this point in 
the corresponding singular variety. When $NS=\calo(1)^n$,
no such normal form can be obtained. In fact, the birational
type of the manifold can be reconstructed from a
complex analytic (and even formal) neighbourhood of $S$.
This was known already to Hartshorne 
(\cite[Theorem 6.7]{_Hartshorne:cohomo-1968_}).

However, one can associate with a quasiline an infinite-dimensional
bundle $\hat\calo_S(M)$ over $S\cong \C P^1$, 
with 
\[ \hat \calo_S(M)=\lim_{\leftarrow} \calo_M/I^n_S\calo_M,
\]
where $I_S$ is an ideal sheaf of $S$. This bundle is called
{\bf a formal completion of $M$ at $S$}.

It is not hard to observe that the space of sections
$H^0(S, \hat \calo_S(M))$ is finite-dimensional, and, moreover,
$H^0(S, \hat \calo_S(M)\otimes_{\calo_S}\calo(i))$
is finite-dimensional for each $i$ 
(\ref{_secti_finite_dim_around_ample_Proposition_}).
For a long time, we expected that the algebra
$A_S:=\bigoplus_i H^0(S, \hat \calo_S(M)\otimes_{\calo_S}\calo(i))$
would provide a sort of an algebraic ``normal form''
of $S$ in $M$, in such a way that the algebraic
structure on $\hat \calo_S(M)$ can be reconstructed
from this ring.

In this paper, we show that this approach works
when $M$ is a Moishezon manifold, and $M$ can be
reconstructed from the ring $A_S$, up to birational
isomorphism (Subsection \ref{_Moishezon_embedded_Subsection_};
this is not very surprising due to the
above-mentioned theorem of Hartshorne, 
\cite[Theorem 6.7]{_Hartshorne:cohomo-1968_}).

For non-Moishezon $M$, such as a twistor space
of a simply connected, compact hyperk\"ahler manifold,
this conjecture is spectacularly wrong. In this case,
$A_S= \bigoplus_i H^0(\C P^1, \calo(i))$, and this ring
has no information about $M$ whatsoever 
(\ref{_algebra_dime_compa_twi_Proposition_}).

\subsection{Holography principle}

The main technical tool of the present paper is the 
following theorem, called {\bf the holography principle}.
Recall that {\bf an ample rational curve}
on a complex variety $M$ is a smooth curve $S\cong \C P^1\subset M$
such that the normal bundle $NS$ is decomposed as $NS=\bigoplus_k\calo(i_k)$,
with all $i_k>0$. 

\hfill

\theorem\label{_hologra_intro_Theorem_}
Let $S\subset M$ be an ample curve in 
a simply connected complex manifold,
and $U$ its connected neighbourhood. Suppose that $M$ is
covered by smooth complex-analytic deformations of $S$.
Then
\begin{description}
\item[(i)] For any holomorphic vector bundle
$B$ on $M$, the restriction map \[ H^0(M,B) \arrow H^0(U,B)\]
is an isomorphism.
\item[(ii)] Let $\Mer(M)$, $\Mer(U)$ be the fields of
meromorphic functions on $M$ and $U$. Then the restriction
map $\Mer(M)\arrow \Mer(U)$ is an isomorphism.
\end{description}

{\bf Proof:} For line bundles, 
\ref{_hologra_intro_Theorem_} (i)
is implied by \ref{_holo_line_Theorem_}.
For general vector bundles, an elegant
argument is given by D. Kaledin in the 
appendix to this paper (Section \ref{_section_Appendix_}).
Holography principle for meromorphic functions is 
proven in \ref{_holo_for_mero_functions_Theorem_}. \endproof

\hfill

The holography principle is not very surprising of
one looks at the neighbourhood of $S$ from the point of
view of complex analysis. The normal bundle to $S$
is obviously positive. Choose a Hermitian metric
$h$ on a neighbourhood of $S$ such that the Chern
connection on $NS$ induced by $h$ has positive curvature.
Let $d_S:\; M \arrow \R$ be the Riemannian distance to $S$
in this Hermitian metric. Since $d_S$ around $S$
is close to the distance in $NS$, the form
$dd^c d_S$ has $n-1$ positive and 1 negative 
eigenvalue in a sufficiently small neighbourhood
of $S$ (\cite{_Demailly:pseudoconvex-concave_}).
This means that $S$ has a neighbourhood $U$ with 
a smooth boundary $\6U$ such that the Levi form
on $\6U$ has one negative and $\dim M-2$
positive eigenvalues. Then, a holomorphic function
defined on an open subset of $U$ and 
continuous on $\6 U$ can be extended outside 
of a boundary; at least, this is the expectation
one has from the solution of the Levi problem.

\subsection{Moishezon twistor spaces}

\definition
Let $M$ be a compact complex manifold.
Define {\bf the algebraic dimension} as $a(M):= \deg_{tr}\Mer(M)$,
where $\deg_{tr}\Mer(M)$ denotes the transcendence degree of
the field of global meromorphic functions on $M$.

\hfill

\definition
A {\bf Moishezon variety} is a compact 
complex variety satisfying $a(M)=\dim M$.

\hfill

The notion of a Moishezon manifold, as it is usually stated,
makes no sense for non-compact varieties. Indeed,
$\deg_{tr}\Mer(M)=\infty$ even if $M$ is an open disk.

However, when $M$ contains an ample curve, the situation 
changes drastically.

\hfill

\theorem\label{_twi_Moishe_intro_Theorem_}
Let $M$ be a complex manifold containing an ample 
rational curve. Then $a(M)\leq \dim M$. Moreover,
if $a(M)=\dim M$, there exists an open embedding
of $M$ to a Moishezon variety $M'$ which satisfies
$\Mer(M')=\Mer(M)$.

\hfill

{\bf Proof:} See \ref{_ample_curve_embe_Moish_Theorem_}.
\endproof

\hfill

Let now $M$ be a twistor space of $X$.
Here, $X$ can be either a hypercomplex
(hyperk\"ahler) manifold, a quaternionic, or 
quaternionic-K\"ahler manifold, or Riemannian
anti-selfdual manifold. We are not very specific,
because the only thing about $M$ which is used
is existence of a large number of quasilines.
The twistor spaces are complex manifold covered
by quasilines, usually non-K\"ahler and
non-quasiprojective (see \ref{_algebra_dime_compa_twi_Proposition_}
 and \ref{_twi_non_quasiproje_Theorem_}).

Since $M=\Tw(X)$ is covered by quasilines
which are by definition ample, we can apply 
\ref{_Harshorne_formal_Theorem_}, and obtain that
$\deg_{tr}\Mer(M)\leq \dim M$.
We call $M$ {\bf a Moishezon twistor space}
if $\deg_{tr}\Mer(M)= \dim M$. This is equivalent
to an existence of an open embedding $M\arrow M'$
of $M$ to a Moishezon manifold (\ref{_twi_Moishe_intro_Theorem_}). 

Moishezon twistor spaces for compact 4-dimensional
anti-selfdual manifolds were discovered by Y.-S. Poon in 
\cite{_Poon:twistors_}, and much studied since then.
Structure theorems about such manifolds were obtained
by F. Campana, \cite{_Campana:twistor_}, and partial
classification results by N. Honda (see e.g. 
\cite{_Honda:classification_}).

The definition given above extends the class
of ``Moishezon twistor manifolds'' significantly.
Let, for instance, $M:=V\3 G$ be a hyperk\"ahler manifold 
which can be obtained using the hyperk\"ahler reduction,
where $V$ is a flat hyperk\"ahler manifold, and $G$ a compact
Lie group (such as the Nakajima quiver variety). Then
the twistor space $\Tw(M)$ is always Moishezon 
(\ref{_quotients_Moish_Theorem_}).

\bigskip

{\bf Acknowledgements:}
Many thanks to Hans-Joachim Hein, Claude LeBrun, Nobuhiro Honda 
and Dima Ka\-le\-din for
interesting discussions on the subject of this article.


\section{Twistor spaces and ample rational curves}



\subsection{Hyperk\"ahler manifolds}


\definition\label{_hyperco_Definition_}
Let $M$ be a manifold, and $I,J,K\in \End(TM)$
endomorphisms of the tangent bundle satisfying the
quaternionic relation
\[
I^2=J^2=K^2=IJK=-\Id_{TM}.
\]
The manifold $(M,I,J,K)$ is called {\bf hypercomplex}
if the almost complex structures $I$, $J$, $K$
are integrable. If, in addition, $M$
is equipped with a Riemannian metric $g$ which
is K\"ahler with respect to $I,J,K$, the
manifold $(M,I,J,K,g)$
is called {\bf hyperk\"ahler}. 

Consider the K\"ahler forms $\omega_I, \omega_J, \omega_K$
on $M$:
\[
\omega_I(\cdot, \cdot):= g(\cdot, I\cdot), \ \
\omega_J(\cdot, \cdot):= g(\cdot, J\cdot), \ \
\omega_K(\cdot, \cdot):= g(\cdot, K\cdot).
\]
An elementary linear-algebraic calculation implies
that the 2-form $\Omega:=\omega_J+\1\omega_K$ is of Hodge type $(2,0)$
on $(M,I)$. This form is clearly closed and
non-degenerate, hence it is a holomorphic symplectic form.

In algebraic geometry, the word ``hyperk\"ahler''
is essentially synonymous with ``holomorphically
symplectic'', due to the following theorem, which is
implied by Yau's solution of Calabi conjecture.

\hfill

\theorem\label{_Calabi-Yau_Theorem_}
Let $(M,I)$ be a compact, K\"ahler, holomorphically
symplectic manifold. Then there exists a unique
hyperk\"ahler metric on $(M,I)$ with the same K\"ahler class.

{\bf Proof:} See \cite{_Yau:Calabi-Yau_}, \cite{_Besse:Einst_Manifo_}.
\endproof

\hfill

\remark
The hyperk\"ahler metric is unique, but there could
be several hyperk\"ahler structures compatible with
a given hyperk\"ahler metric on $(M,I)$, if the holonomy
of its Levi-Civita connection is strictly less than $\Sp(n)$.

\hfill

\definition \label{_indu_comple_str_Definition_} 
Let $M$ be a hypercomplex manifold, and $L$ a quaternion satisfying
$L^2=-1$. Then $L=aI+bJ+cK$, $a^2+b^2+c^2=1$.
The corresponding complex structure 
on $M$ is called {\bf an induced complex structure}. 
The space $M$, considered as a complex manifold, is denoted by $(M, L)$. 
The set of induced complex structures is naturally identified
with $S^2$, which we often consider as $\C P^1$ with the
standard complex structure.

\hfill

\definition\label{_trianalytic_Definition_} 
(\cite{_Verbitsky:trianalyt_})
Let $X\subset M$ be a closed subset of a hyperk\"ahler
manifold $M$. Then $X$ is
called {\bf trianalytic} if $X$ is a complex analytic subset 
of $(M,L)$ for every induced complex structure $L$.
 
\hfill

Trianalytic subvarieties were a subject of a 
long study. Most importantly, consider a generic
induced complex structure $L$ on $M$. Then all closed complex
subvarieties of $(M, L)$ are trianalytic. 
Moreover,  a trianalytic subvariety can be canonically
desingularized (\cite{_Verbitsky:hypercomple_}),
and this desingularization is hyperk\"ahler.

\hfill

\theorem\label{_triana_generic_Theorem_}
(\cite{_Verbitsky:trianalyt_}, \cite{_Verbitsky:subvar_})
Let $(M,I,J,K)$ be a hyperk\"ahler manifold (not necessarily
compact). Then there exists a countable subset $R\subset \C P^1$,
such that for any induced complex structure $L\notin R$,
all compact complex subvarieties of $(M,L)$ are
trianalytic.
\endproof

\hfill

\remark
For hypercomplex manifolds, \ref{_triana_generic_Theorem_}
is (generally speaking) false, though for manifolds
with trivial canonical bundle a weaker form of
this result was obtained (\cite{_Soldatenkov_Verbitsky:triana_}).


\subsection{Twistor spaces}


\definition
Let $M$ be a Riemannian 4-manifold. Consider the action of
the Hodge $*$-operator: $*:\; \Lambda^2 M \arrow \Lambda^2 M$.
Since $*^2 =1$, the eigenvalues are $\pm 1$, and one
has a decomposition $\Lambda^2 M = \Lambda^+ M \oplus \Lambda^- M$
onto {\bf selfdual} ($*\eta=\eta$) and 
{\bf anti-selfdual} ($*\eta=-\eta$) forms.

\hfill

\remark If one changes the orientation of $M$,
leaving  metric the same, $\Lambda^+ M$ and $\Lambda^- M$
are exchanged. Therefore, their dimensions
are equal, and $\dim \Lambda^2 M=6$ 
implies $\dim \Lambda^\pm(M)=3$.

\hfill

\remark
Using the isomorphism $\Lambda^2 M = \goth{so} (TM)$,
we interpret $\eta \in \Lambda^2_m M$ as an endomorphisms of $T_mM$.
Then the unit vectors $\eta \in \Lambda^+_mM$ correspond
to oriented, orthogonal complex structures on $T_m M$.

\hfill

\definition
Let $\Tw(M):=S\Lambda^+ M$ be the set of unit vectors in $\Lambda^+M$.
At each point $(m,s)\in \Tw(M)$, consider the decomposition
$T_{m,s}\Tw(M)= T_m M \oplus T_sS\Lambda^+_m M$, induced by the
Levi-Civita connection.
Let $I_s$ be the complex structure on $T_mM$ induced by $s$,
$I_{S\Lambda^+_m M}$ the complex structure on $S\Lambda^+_m M=S^2$
induced by the metric and orientation, and 
\[ {\cal I}:\; T_{m,s}\Tw(M)\arrow T_{m,s}\Tw(M)\]
be equal to ${\cal I}_s \oplus I_{S\Lambda^+_m M}$.
An almost complex manifold $(\Tw(M), {\cal I})$
is called {\bf the twistor space} of $M$.

\hfill

Given a hyperk\"ahler or hypercomplex manifold (\ref{_hyperco_Definition_}),
one defines its twistor space in a similar manner.

\hfill

\definition
A {\bf twistor space} $\Tw(M)$ of a hypercomplex
manifold $M$ is  $S^2 \times M$ equipped
with a complex structure which is defined
as follows.  Consider the complex structure $I_m:T_mM \to T_mM$ 
on $M$ induced by $J \in S^2 \subset {\Bbb H}$. Let $I_J$
denote the complex structure on $S^2 = \C P^1$.
The operator $I_{\Tw} = I_m \oplus I_J:T_x\Tw(M) \to T_x\Tw(M)$ 
satisfies $I_{\Tw} ^2 = -\Id$.  It defines 
an almost complex structure on $\Tw(M)$.

\hfill

The almost complex structure on 
the twistor space of a Riemannian 4-manifold $X$
is integrable whenever $X$ is anti-selfdual
(\cite{_Atiyah_Hitchin_Singer_})
For a hypercomplex manifold it is integrable as well
(\cite{_Kaledin:twistor_}).
Twistor spaces are the main example of the
geometries we are working with.


\subsection{Geometry of twistor spaces}


\proposition\label{_algebra_dime_compa_twi_Proposition_}
Let $\Tw(M)$ be a twistor space of a compact hyperk\"ahler
manifold. Then 
\begin{description}
\item[(i)] $\Tw(M)$ is non-K\"ahler.
\item[(ii)] The algebraic dimension of $\Tw(M)$ is 1.
\end{description}

{\bf Proof of (i):}
Let $\omega$ be the standard Hermitian form of $\Tw(M)$.
Then $dd^c\omega$ is a positive (2,2)-form (\cite[(8.2)]{_NHYM_}).
For any K\"ahler form
$\omega_0$ on $\Tw(M)$, this would imply 
\[ 
 \int_{\Tw(M)}d\left(\omega_0^{\dim_\C M-1} \wedge d^c\omega\right)=
 \int_{\Tw(M)}\omega_0^{\dim_\C M-1} \wedge dd^c\omega>0,
\]
which is impossible by Stokes' theorem.

\hfill

{\bf Proof of (ii):}
See \ref{_dim_tr_twistor_compact_Theorem_}.
\endproof


\subsection{Quasilines and ample curves in the twistor spaces}


\definition
{\bf An ample rational curve} 
on a complex manifold $M$ is
a smooth curve $S \cong \C P^1\subset M$ such that 
its normal bundle $N\!S$ satisfies 
$N\!S=\bigoplus_{k=1}^{n-1}\calo(i_k)$, with all $i_k >0$
(see \cite{_Kollar:curves_}).
It is called {\bf a quasiline} if all
$i_k=1$.

\hfill

\claim
Let $M$ be a twistor
space of a hyperk\"ahler or 4-dimensional ASD manifold, 
$m\in M$ a point, and $S_m$
the corresponding $S^2$ in $\Tw(M)=\C P^1 \times M$. Then $S_m$
  is a quasiline.

\hfill

{\bf Proof:}
Since the claim is essentially
infinitesimal, it suffices to check it when $M$ is flat.
Then $\Tw(M)= \Tot (\calo(1)^{\oplus 2n})
\cong \C P^{2n+1} \backslash \C P^{2n-1}$, and
$S_m$ is a section of $\calo(1)^{\oplus 2n}$.
\endproof

\hfill

Existence of quasilines in twistor spaces is a very strong
condition, and can be used to obtain all kinds of
geometric information; for example, see \cite{_Campana:twistor_} and
\cite{_Verbitsky:rational_twi_}.


\section[Holography principle for manifolds with ample 
rational curves]{Holography principle for manifolds \\  with ample 
rational curves}


\subsection{Holography principle for line bundles}

Throughout this paper, all 
neighbourhoods and manifolds are silently assumed
to be connected. One of the main results of the present paper
is the following theorem. 

\hfill

\theorem\label{_holo_line_Theorem_}
(holography principle for line bundles)
Let $S\subset M$ be an ample
rational curve in a simply connected complex manifold,
which is covered by smooth, ample deformations of $S$, and
$L$ a holomorphic line bundle on $M$.
Consider an open neighbourhood $U\supset S$.
Then the 
restriction map $H^0(M,L) \arrow H^0(U,L)$ is an isomorphism.

\hfill

We deduce \ref{_holo_line_Theorem_}
from the following local result
(\ref{_holo_line_local_Proposition_}).

\hfill

\remark\label{_ample_defo_neighbo_interse_Remark_}
Since $S$ is an ample curve, $S$ can be deformed
in any normal direction. Therefore, there exists
an open neighbourhood $U\supset S$ which is contained
in a union of the set ${\goth S}$ of all smooth, 
ample deformations of $S$ intersecting $S$.
Further on, we choose this neighbourhood in such a 
manner that any $S_1\in {\goth S}$ can be connected to $S$
by a continuous family of deformations intersecting $S$.

\hfill

\proposition\label{_holo_line_local_Proposition_}
Let $S\subset M$ be an ample
rational curve in a connected complex manifold,
which is covered by deformations of $S$, and
$L$ a holomorphic line bundle on $M$.
Consider a neighbourhood $U\supset S$ which is is contained
in a union of the set of all deformations of $S$
intersecting $S$ (\ref{_ample_defo_neighbo_interse_Remark_}).
Then for any smaller 
open neighbourhood $V\subset U$ of $S$,
the  restriction map $H^0(U,L) \arrow H^0(V,L)$ 
is an isomorphism.

\hfill

{\bf Proof of an implication
``\ref{_holo_line_local_Proposition_} $\Rightarrow$
\ref{_holo_line_Theorem_}''.}

{\bf Step 1:} Choose a 
continuous, connected family $S_b$ of ample curves parametrized
by $B$ such that $\bigcup_{b\in B} S_b = M$, and
choose a tubular neighbourhood $U_b$ for each $S_b$,
continuously depending on $b$. Then the 
intersection $U_b\cap U_{b'}$
for sufficiently close $b,b'$ always contains $S_b$ and
$S_{b'}$. 
By \ref{_holo_line_local_Proposition_}, $U_b$
can be chosen in such a way that
$H^0(U_b\cap U_{b'},L)=H^0(U_b,L)=H^0( U_{b'},L)$.

{\bf Step 2:} Since $B$ is connected, all 
the spaces $H^0(U_b,L)$ are isomorphic, and these isomorphisms
are compatible with the restrictions to the intersections 
$U_b\cap U_{b'}$. Let now $f\in H^0(U_b,L)$, and let
$\tilde M_f$ be the domain of holomorphy for $f$,
that is, a maximal domain (non-ramified over $M$) such that
$f$ admits a holomorphic extension to $\tilde M_f$. 
Since $\cup U_b=M$, and $f$ can be holomorphically
extended to any $U_b$, the domain $\tilde M_f$
is a covering of $M$. Now, \ref{_holo_line_Theorem_}
follows, because $M$ is simply connected.

\subsection{The local holography principle}

To prove \ref{_holo_line_Theorem_}
it remains to prove \ref{_holo_line_local_Proposition_}.
If $\deg L\restrict S<0$, the statement of \ref{_holo_line_local_Proposition_}
is vacuous. Therefore, we may always assume that $\deg L\restrict S\geq 0$,
hence $l\restrict S$ is generated by global sections. 
Let $S(M)$ be the space of deformations of the ample curve
$S$ which remain smooth, and $SS(M)$ the space of pairs
$\{(x, S_1):\ \ S_1\in S(M), x\in S_1\}$. Consider the
natural forgetful maps $\tau_1:\; SS(M)\arrow M$, 
$\tau_2:\; SS(M)\arrow S(M)$, and let $E$ be the bundle
${\tau_2}_*\tau_1^*L$ on $S(M)$. 
Denote by $\deg_SL$ the degree of the restriction of $L$ to
$S$. Since $\dim H^0(S_1, L)=\deg L+1$,
$E$ is a $(\deg L+1)$-dimensional vector bundle.
Given a section $f$ of $L$ on $M$, denote the 
corresponding section of $E$ by $\tilde f:={\tau_2}_*\tau_1^*f$.
When $\deg_S L=d$, the value of $\tilde f$ at $S_1$
is uniquely determined by the restriction of $f$
to any $d+1$ distinct points of $S_1$. Indeed,
$E\restrict{S_1}$ is $d+1$-dimensional, and 
any section $h\in H^0(S_1,L)$ is uniquely determined
by its values in $d+1$ points. This gives a map
\begin{equation}\label{_reconstru_section_Equation_}
L\restrict{z_1}\times L\restrict{z_2}\times ...\times
L\restrict{z_{d+1}} \tilde\arrow H^0(S_1, L).
\end{equation}

Now, let $f$ be a section of $L$ on $V$.
For any $S_1\in S(M)$ intersecting $V$, we choose
$d_1$ distinct points $z_1, ..., z_{d+1}\in S_1 \cap V$,
and consider the section 
$\phi\left(f\restrict{z_1}, f\restrict{z_2}, ...,
f\restrict{z_{d+1}}\right)\in H^0(S_1,L)$
defined using \eqref{_reconstru_section_Equation_}.
When $S_1\subset V$, this section is independent
from the choice of $z_1, ..., z_{d+1}\in S_1$.
Let $R_V$ be a connected component of the set of all
$S_1\in S(M)$ intersecting $V$ and containing $S$.
Since the map 
\[ (z_1,..., z_{d+1})\arrow \phi\left(f\restrict{z_1},
f\restrict{z_2}, ..., f\restrict{z_{d+1}}\right)
\]
is holomorphic and independent from the choice
of $z_1,..., z_{d+1}$ on an open subset of $R_V$,
it is independent of $z_1,..., z_{d+1}$ everywhere on $R_V$.
This gives a section $\tilde f\in H^0(R_V,E)$
extending the section 
$\tilde f:={\tau_2}_*\tau_1^*f\in H^0(S(V),E)$.

By construction, $U$ is contained in a connected
part $U_1$ of the union of all deformations of $S$
intersecting $V$. To extend $f$ from $V$ to $U$,
we use $\tilde f$ to obtain an extension of $f$
to $U_1$, as follows.

Any section $g\in H^0(R_V, E)$ 
gives a function $\psi_g$ mapping a pair $(x,S_1)$,
$x\in S_1\in R_V$
to $g(S_1)\restrict x\in L\restrict x$.
For the section $\tilde f$ constructed above,
$\psi_{\tilde f}(x,S_1)$ is independent from the choice of $S_1$
whenever $S_1$ lies in $V$. The same analytic continuation argument
as above implies that $\psi_{\tilde f}(x,S_1)$ is
independent of $S_1$ everywhere. For any $x\in U_1$,
the set $f\in L\restrict x$ equal to 
$\psi_{\tilde  f}(x,S_1)$, where 
$S_1\in R_V$ is an arbitrary curve passing through $x$.
This gives an extension of $f$ to $U_1$.
\ref{_holo_line_local_Proposition_} is proven.
We finished the proof of \ref{_holo_line_Theorem_}.
\endproof

\subsection{Holography principle for meromorphic functions}

The following theorem is proven in the same way as
\ref{_holo_line_Theorem_}. Given a complex variety $M$,
we denote the field of meromorphic functions on $M$ by $\Mer(M)$.

\hfill

\theorem\label{_holo_for_mero_functions_Theorem_}
Let $S\subset M$ be an ample curve in a 
simply connected complex manifold, and $U\supset S$ a connected
neighbourhood of $S$. Suppose that $M$ is covered by the
union of all smooth, ample deformations of $S$.
Then the restriction map $\Mer(M)\arrow \Mer(U)$ is an isomorphism.

\hfill

{\bf Proof:} The same argument as used to deduce
\ref{_holo_line_Theorem_} from \ref{_holo_line_local_Proposition_}
can be used to reduce
\ref{_holo_for_mero_functions_Theorem_}
to the following statement.

\hfill

\proposition\label{_holo_mero_local_Proposition_}
Let $S\subset M$ be an ample
rational curve in a connected complex manifold,
which is covered by smooth, ample deformations of $S$.
Consider a neighbourhood $U\supset S$ which is contained
in a union ${\goth S}$ of all deformations of $S$
intersecting $S$ (\ref{_ample_defo_neighbo_interse_Remark_}).
Then for any smaller 
open neighbourhood $V\subset U$ of $S$,
the restriction map $\Mer(U) \arrow \Mer(V)$ 
is an isomorphism.

\hfill

{\bf Proof:} Define {\bf the degree} $\deg_S(f)$ of a meromorphic
function $f\in \Mer(V)$ as the degree of the pole divisor
of $f\restrict {S_1}$ for any deformation $S_1$ of $S$ transversal
to the pole divisor of $f$. Denote by $\Mer_d(V)$
the space  meromorphic functions of degree $\leq d$. To prove
\ref{_holo_mero_local_Proposition_} it would suffice to
show that the restriction map $\Mer_d(U) \arrow \Mer_d(V)$ 
is an isomorphism, for all $d$.

For each rational curve $S_1$, a 
degree $\leq d$ meromorphic function is uniquely determined
by its values in any $d+1$ distinct points on $S_1$.
Given a meromorphic function $f\in \Mer_d(V)$,
and a deformation $S_1$ of $S$ intersecting $V$, we can 
extend $f\restrict {S_1\cap V}$ to a degree $d$
meromorphic function $f_1$ on $S_1$ by computing its values
at $d+1$ distinct points
$z_1, ..., z_{d+1}$ of $S_1\cap V$. Whenever $S_1$ is in $V$, this 
procedure gives $f\restrict {S_1}$. By analytic continuation,
the values of $f_1(z)$ at any $z\in S_1$ are independent
from the choice of $z_i$ and $S_1$. 

We have shown that $f_1$ is a well-defined meromorphic
function on the union $U_1$ of all deformations of $S$ 
intersecting $V$. By construction, on $V$ we have $f_1=f$.
Since $U_1$ contains $U$, this implies
that $f$ can be extended from $V$ to $U$.
\endproof


\section{Moishezon twistor spaces}


\subsection{Vector bundles in a
  neighbourhood of an ample curve}

\proposition\label{_secti_finite_dim_around_ample_Proposition_}
Let $S\subset M$ be an ample rational curve,
$U\supset S$ its  neighbourhood, and
$B$ a holomorphic bundle on $U$. Then $H^0(U,B)$
is finite-dimensional.

\hfill

{\bf Proof. Step 1:}
Let $I_S\subset \calo_U$ be the ideal sheaf of $S$, and 
$J^r_S(B):= B / I^{r+1}_SB$ the sheaf of $r$-jets of the sections
of $B$.  Since $I^r_S/I^{r+1}_S=\Sym^r(N^*S)$,
the bundle $I^r_S/I^{r+1}_S$ is a direct sum of $\calo(k_i)$
with $k_i < -r$. Therefore, $H^0(I^r_SB/I^{r+1}_SB)=0$
for $r$ sufficiently big.

\hfill

{\bf Step 2:}
 The sheaf $I^{r}_SB$ admits a filtration
$I^{r}_SB\supset I^{r+1}_SB\supset I^{r+2}_SB\supset ...$
with associated graded sheaves $I^r_SB/I^{r+1}_SB$ 
having no sections for $r\gg 0$. Therefore, $H^0(I^{r}_SB)=0$
for sufficiently byg $r$.

\hfill

{\bf Step 3:}
If $H^0(U,B)$ is infinite-dimensional, the 
map \[ H^0(U,B)\arrow H^0(B/I^r_S(B))\] cannot be injective.
Then, for each $r$, there exists a non-zero section
with vanishing $r$-jet: $f_r\in H^0(U,I^{r+1}_SB)$.
This is impossible, as shown in Step 2. \endproof

\hfill

This proof is effective, and gives the following
bound on the dimension of the space of sections of $B$.

\hfill

\corollary\label{_space_of_sec_effective_Corollary_}
Let $S\subset M$ be an ample rational curve,
$U\supset S$ its  neighbourhood, and
$B$ a holomorphic bundle on $U$. Then 
\[ \dim H^0(U,B)\leq \dim H^0\left(S,
   \bigoplus_d \Sym^d(N^*S)\otimes_{\calo_S}B\restrict S\right).
\]
\endproof

\hfill

The same argument can be applied to degree $d$ meromorphic
functions. Recall that {\bf the degree} $\deg_S(f)$ of a meromorphic
function $f\in \Mer(V)$ is the degree of the pole divisor
of $f\restrict {S_1}$ for any deformation of $S$ transversal
to the pole divisor of $f$. We denote the space of
meromorphic functions of degree $\leq d$ on $U$ by $\Mer_d(U)$.

\hfill

\corollary\label{_space_of_mer_effective_Corollary_}
Let $S\subset M$ be an ample rational curve, and
$U\supset S$ its neighbourhood. Then 
\begin{equation}
\label{_mero_dime_bound_Equation_}
\dim \Mer_d(U)\leq \dim H^0\left(S,
   \Sym^{\leq d}(N^*S)\otimes_{\calo_S}\calo(d)\right).
\end{equation}
\endproof

\subsection{Algebraic dimension of the field of meromorphic functions}

The conormal bundle $N^*\!S$ is negative; clearly, the dimension 
\eqref{_mero_dime_bound_Equation_} is maximal when $S$ is
a quasiline, and $N^*\!S=\oplus \calo(1)$. In this case, the
bound \eqref{_mero_dime_bound_Equation_} is realized for
a rational line in $\C P^n$. Indeed, for a 
rational line $S$ in $\C P^n$, the sheaf of algebraic functions
in a neighbourhood of $S$ is isomorphic to 
of $\bigoplus_i \Sym^i(N^*\!S)$. This implies the following
simple numerical result.

\hfill

\claim
Let $S\subset M$ be an ample rational curve, and
$U\supset S$ its neighbourhood. Then 
\begin{equation}
\label{_mero_dime_bound_expli_Equation_}
\dim \Mer_d(U)\leq \dim H^0(\C P^n,\calo(d)).
\end{equation}
where $n=\dim M$.
\endproof

\hfill

\corollary\label{_deg_tr_twistor_Theorem_}
Let $M$ be a complex variety containing an
ample rational curve. Then the transcendence degree of
$\Mer(M)$ satisfies $\deg_{tr}\Mer(M) \leq \dim M$.

\hfill

{\bf Proof:} Consider the graded 
ring $\bigoplus_d \Mer_d(M)$. Since
$\dim \Mer_d(U)\leq \dim H^0(\C P^n,\calo(d))$,
the Krull dimension of this ring is $\leq n$.
Therefore, the transcendence degree of its
ring of fractions is also bounded by $n$.
\endproof

\hfill

This observation is not new: it was known already 
to Hartshorne (in the context of formal neighbourhoods).
Applied to complex analytic spaces, Harts\-horne's theorem
can be stated as follows.

\hfill

\theorem \label{_Harshorne_formal_Theorem_}
Let $S\subset M$ be a connected, positive-dimensional,
smooth subvariety in a complex manifold.
Assume that the normal bundle of $S$ is ample. 
Then the transcendence degree of the field $\Mer(M)$
of meromorphic functions is no bigger
than the dimension of $M$:
\[ \deg_{tr}\Mer(M)\leq \dim(M).
\]
Moreover, if equality is reached, $\Mer(M)$ is a finitely
generated extension of $\C$.

\hfill

{\bf Proof:} \cite[Theorem 6.7]{_Hartshorne:cohomo-1968_};
see also \cite{_Kebekus_Sola_Toma_}. \endproof


\subsection{Moishezon manifolds and ample rational curves}
\label{_Moishezon_embedded_Subsection_}


Hartshorne's methods are already sufficient to prove
the following general result.

\hfill

\proposition\label{_ample_curve_local_embe_Proposition_}
Let $S\subset M$ be an ample rational 
curve in a simply connected complex manifold.
Assume that $\deg_{tr}\Mer(M)= \dim(M).$
Then there exists a meromorphic map
$\phi:\; M\arrow Z_0$ to an open subset of 
a projective variety $M$, which is 
bijective onto its image outside of a complex
analytic subset of positive codimension.

\hfill

{\bf Proof:} By Hartshorne's theorem 
(\ref{_Harshorne_formal_Theorem_}),
$\Mer(M)$ is a finitely
generated extension of $\C$. Let $\xi_1, ..., \xi_N$ be generators 
of $\Mer(M)$, $D_i$ their pole divisors, and $L_i:=\calo(D_i)$
the corresponding line bundles. Then $\xi_i$ can be considered as
sections of $L_i$, and $\xi_1, ..., \xi_N$ -- as sections
of $L:=\bigotimes_i L_i$. Consider now the subring of
$\bigoplus_d H^0(M, L^d)$ generated by $\xi_i$, and 
let $Z$ be its spectre. Clearly, $\dim M=\dim Z$, $Z$ is projective,
and the natural rational map
$M\stackrel \phi\arrow Z$ induces an isomorphism $\Mer(Z)\tilde \arrow \Mer(M)$.

Let $\tilde M$ be a resolution of the base set of $\phi$,
such that $\tilde M\stackrel {\tilde \phi}\arrow Z$ is holomorphic.
If $\tilde \phi$ is ramified at some divisor $D$ in $\phi(\tilde M)$, 
this divisor can be extended to $Z$ using 
\ref{_holo_for_mero_functions_Theorem_} applied to 
$U=\phi(\tilde M)\subset Z$. Taking the corresponding ramified 
covering $\tilde Z$ of $Z$, we obtain a holomorphic map from
$\tilde M$ to $\tilde Z$, which is impossible, because $\Mer(\tilde Z)$
is strictly bigger than $\Mer(Z)$, and $\Mer(\tilde M)=\Mer(Z)$. 
Therefore, $\tilde \phi$ is 
bijective to its impage at its general point. 
Openness of its image is a general property
of bimeromorphic maps. \endproof

\hfill

\theorem\label{_ample_curve_embe_Moish_Theorem_}
Let $S\subset M$ be an ample rational 
curve in a simply connected complex manifold.
Assume that $\deg_{tr}\Mer(M)= \dim(M).$
Then there exists an open embedding of $M$ 
to a Moishezon variety.

\hfill

{\bf Proof:} From \ref{_ample_curve_local_embe_Proposition_}
we obtain a line bundle $L$ on $M$ inducing a bimeromorphic
map $\phi:\; M\arrow Z_0$ to an open subset of 
a projective variety $Z$. Resolving the base points of
the inverse map if necessary, we may assume that
the inverse map $\psi:\; Z_0 \arrow M$ is holomorphic.
Then, $M$ is obtained from $Z_0$ by blowing down 
a certain number of exceptional subvarieties $E_i$, 
obtained as common zero sets of a 
certain number of meromorphic functions.
Applying \ref{_holo_for_mero_functions_Theorem_}
to $Z_0\subset Z$, we extend the meromorphic 
functions and the corresponding subvarieties 
$E_i$ and obtain closed exceptional subvarieties
 $E_i'\subset Z$. Blowing these down, we obtain a Moishezon
variety which contains $M$ as an open subset.
\endproof


\subsection{Moishezon twistor spaces}


\definition\label{_twistor_Moish_Definition_}
Let $M=\Tw(X)$ be a twistor space of a 
simply connected hyperk\"ahler, hypercomplex,
quaternionic or or 4-dimensional anti-selfdual manifold,
not necessarily compact. We say that $M$
is {\bf a Moishezon twistor space} if 
$\deg_{tr}\Mer(M)= \dim(M)$ (see \ref{_deg_tr_twistor_Theorem_}).

\hfill

{}From \ref{_ample_curve_embe_Moish_Theorem_},
we immediately obtain the following corollary.

\hfill

\corollary
Let $M$ be a Moishezon twistor space.
Then $M$ admits an open embedding to a Moishezon variety $M_1$.
Moreover, $M_1$ is unique up to a bimeromorphic equivalence.

\hfill

{\bf Proof:} The open embedding to a 
Moishezon variety follows from 
\ref{_ample_curve_embe_Moish_Theorem_},
and its uniqueness is implied by an isomorphism
$\Mer(M)=\Mer(M_1)$ (\ref{_holo_for_mero_functions_Theorem_}).
\endproof

\hfill

It is easy to construct an example of a twistor
manifold which does not belong to this class.
The twistor space of a K3 surface, and, more generally,
any compact hyperk\"ahler manifold is never Moishezon.

\hfill

\theorem\label{_dim_tr_twistor_compact_Theorem_}
Let $M$ be a compact hyperk\"ahler manifold,
and $\Tw(M)$ its twistor space. Then $\deg_{tr}\Mer(\Tw(M))=1$.

\hfill

{\bf Proof:} Let $Z\subset \Tw(M)$ be any divisor,
and $R\subset \C P^1$ a countable subset constructed
in \ref{_triana_generic_Theorem_}. 
For any induced complex structure $L\notin R$,
all complex subvarieties of $(M,L)$ are even-dimensional. 
Therefore, $Z$ intersects the twistor fiber 
$(M,L)=\pi^{-1}(L)\subset \Tw(M)$
non-transversally, or not at all. By Thom's transversality theorem, 
the intersection $Z$ with all fibers of $\pi$ except
a finite number is transversal. This means that
$Z$ can intersect only finitely many of the fibers
of $\pi$. However, all these fibers are irreducible
divisors. Therefore, $D$ is a union 
of several fibers of $\pi$. Since a meromorphic
function is uniquely determined by its pole or
zero divisor, all meromorphic functions on $\Tw(M)$
are pull-backs of meromorphic functions on $\C P^1$.
\endproof

\hfill

\remark
When $M$ is  hypercomplex, no effective
bounds on the transcendence degree $\deg_{tr}\Mer(\Tw(M))$ are known.
We conjecture (based on empirical evidence)
that the twistor space of a compact
hypercomplex manifold is not Moishezon, but 
this conjecture seems to be difficult.


\subsection{Twistor spaces and hyperk\"ahler reduction}


We recall the definition of hyperk\"ahler reduction,
following \cite{_HKLR_} and  \cite{_Nakajima_}. This material
is fairly standard.

\hfill

We denote the Lie derivative along a vector field
as $\Lie_x:\; \Lambda^i M \arrow \Lambda^i M$,
and contraction with a vector field by 
$i_x:\; \Lambda^i M \arrow \Lambda^{i-1} M$.
Recall the Cartan's formula:
\[
d\circ i_x + i_x \circ d =\Lie_x
\]

Let $(M,\omega)$ be a symplectic manifold, $G$ a 
Lie group acting on $M$ by symplectomorphisms, and $\goth g$
its Lie algebra. For any $g\in {\goth g}$,
denote by $\rho_g$ the corresponding vector field.
Cartan's formula gives $\Lie_{\rho_g}\omega=0$, hence
$d(i_{\rho_g}(\omega))=0$. We obtain that 
$i_{\rho_g}(\omega)$ is closed, for any $g\in {\goth g}$.

\hfill

\definition
{\bf A Hamiltonian} of $g\in {\goth g}$ is a function
$h$ on $M$ such that $dh=i_{\rho_g}(\omega)$.

\hfill

\definition
$(M,\omega)$ be a symplectic manifold, $G$ a 
Lie group acting on $M$ by symplectomorphisms.
A moment map $\mu$ of this action is a linear map
${\goth g}\arrow C^\infty M$
associating to each $g\in \g$ its Hamiltonian.

\hfill

\remark
It is more convenient to consider $\mu$
as an element of ${\goth g}^* \otimes_\R C^\infty M$,
or, as it is usually done,  a function on $M$
with values in ${\goth g}^*$. 

\hfill

\remark
Note that the moment map always exists, if $M$ 
is simply connected.

\hfill

\definition 
A moment map $M \arrow  {\goth g}^*$
is called {\bf equivariant}
if it is equivariant with respect to the 
coadjoint action of $G$ on ${\goth g}^*$.

\hfill

\remark
$M\stackrel\mu \arrow  {\goth g}^*$ is a moment map
if and only if for all $g\in {\goth g}$, 
$\langle d\mu,g\rangle= i_{\rho_g}(\omega)$.
Therefore,  a moment map is defined up to 
a constant ${\goth g}^*$-valued function.
An equivariant moment map is is defined up to 
 a constant ${\goth g}^*$-valued function
which is $G$-invariant.

\hfill

\definition A $G$-invariant $c\in \goth g^*$ is called
{\bf central}.

\hfill

\claim
An equivariant moment map exists whenever $H^1(G, {\goth g}^*)=0$.
In particular, if $G$ is reductive and $M$ is simply connected,
an equivariant moment map is always possible to define.
\endproof

\hfill

\definition
Let $(M,\omega)$ be a symplectic manifold, $G$ a 
compact Lie group acting on $M$ by symplectomorphisms,
$M\stackrel\mu \arrow  {\goth g}^*$ an equivariant
moment map, and $c\in {\goth g}^*$ a central element.
The quotient $\mu^{-1}(c)/G$ is called the {\bf 
symplectic reduction} of $M$, denoted by $M\2 G$.

\hfill

\claim
The symplectic quotient $M\2 G$ is a symplectic
manifold of dimension $\dim M - 2 \dim G$.
\endproof

\hfill

\theorem
Let $(M,I, \omega)$ be a K\"ahler manifold,
$G_\C$ a complex reductive Lie group acting on $M$ by
holomorphic automorphisms, and $G$ is a compact form of $G_\C$
acting isometrically.  Then $M\2 G$ is a K\"ahler orbifold.
\endproof

\hfill

\remark
In such a situation, $M\2 G$
is called {\bf the K\"ahler quotient},
or {\bf  GIT quotient}. 

\hfill

\remark
The points of $M\2 G$ are in bijective correspondence
with the orbits of $G_\C$ which intersect
$\mu^{-1}(c)$. Such orbits are called
{\bf polystable}, and the intersection
of a $G_\C$-orbit with $\mu^{-1}(c)$ is a $G$-orbit.

\hfill

\definition
Let $G$ be a compact Lie group, $\rho$ its action 
on a hyperk\"ahler manifold $M$ by hyperk\"ahler isometries, and
$\g^*$ a dual space to its Lie algebra. {\bf A 
hyperk\"ahler moment map} is a $G$-equivariant 
smooth map $\mu: M\to\g^*\otimes\R^3$ such that
 $\langle \mu_i(v),g \rangle = \omega_i(v,d\rho(g))$, 
for every $v\in TM$, $g\in\g$ and $i=1,2,3$,
where $\omega_i$ are three K\"ahler forms associated 
with the hyperk\"ahler structure.

\hfill

\definition
Let $\xi_1,\xi_2,\xi_3$ be three $G$-invariant vectors
in $\g^*$.
The quotient manifold $M\3 G := \mu^{-1}(\xi_1,\xi_2,\xi_3)/G$ 
is called {\bf the hyperk\"ahler quotient} of $M$.

\hfill

\theorem (\cite{_HKLR_}) The quotient $M\3 G$
is hyperkaehler.

\hfill

{\bf Proof:} We sketch the proof of
Hitchin-Karlhede-Lindstr\"om-Ro\v cek theorem, because
we make use of it further on.

Let $\Omega:=\omega_J+ \1\omega_K$. This is a holomorphic
symplectic (2,0)-form on $(M,I)$.
Let $\mu_J, \mu_K$ be the moment map associated with 
$\omega_J, \omega_K$, and $\mu_\C:=  \mu_J+ \1\mu_K$.
Then $\langle d\mu_\C,g\rangle= i_{\rho_g}(\Omega)$
Therefore, $d\mu_\C\in \Lambda^{1,0}(M,I)\otimes {\goth g}^*$.
This implies that the map $\mu_\C$ is
holomorphic. It is called {\bf a holomorphic moment map.}

By definition, $M\3 G=\mu_{\C}^{-1}(c)\2 G$, 
where $c\in {\goth g}^*\otimes_\R \C$
is a central element.
This is a K\"ahler manifold, because it is a K\"ahler
quotient of a K\"ahler manifold.

 We obtain 3 complex structures
$I,J,K$ on the hyperk\"ahler quotient $M\3 G$.
They are compatible in the usual way, as seen from a 
simple local computation.
\endproof

\hfill

\theorem\label{_quotients_Moish_Theorem_}
Let $V$ be a quaterionic Hermitian vector space, and
$G\subset \Sp(V)$ a compact Lie group acting on $V$ by quaternionic
isometries. Denote by $M$ the hyperk\"ahler reduction of $V$.
Then $\Tw(M)$ is a Moishezon twistor space, in the sense
of \ref{_twistor_Moish_Definition_}. 

\hfill

{\bf Proof:}
The holomorphic symplectic form on $(V,I)$ depends on $I$
holomorphically, giving a section 
\[ \Omega_{tw}\in\Omega^2_\pi(\Tw(V))\otimes_{\calo_{\Tw(V)}}\pi^*\calo(2).\]
Here, $\pi:\; \Tw(V)\arrow \C P^1$ is the twistor projection, and
$\Omega^2_\pi(\Tw(V))$ the sheaf of fiberwise holomorphic 2-forms.
Consider the fiberwise holomorphic moment map given by this form,
$\mu_{tw}:\; \Tw(V)\arrow \Tot({\goth g}^*\otimes_\C \calo(2))$.
Replacing the moment map by its translate, we can always assume that
$M\3 G=\mu_{\C}^{-1}(0)\2 G$ (that is, we assume that 
the central vector $c$ used to define $M\3 G$ vanishes).
Then $\Tw(M)$ is obtained
as the space of polystable $G_\C$-orbits in
$\mu_{tw}^{-1}(0)\subset \Tw(V)$. Here, ``polystability''
of an orbit is understood as the non-emptiness of the intersection
of this orbit with $\mu^{-1}(0)$; the set of such orbits is open.

The space 
$\Tw(V)=\C P^{2n+1}\backslash \C P^{2n-1}$
is quasiprojective. Averaging the ring of rational functions with
$G$, we obtain that the field $G_\C$-invariant rational
functions on $\Tw(V)$ has dimension 
$\dim \mu_{tw}^{-1}(0)-\dim G_\C=\dim \Tw(M)$,
and hence $\Tw(M)$ is Moishezon. \endproof

\hfill

\corollary
Let $U$ be an open subset of a compact, simply
connected hyperk\"ahler manifold, and $U'$
an open subset of a hyperk\"ahler manifold obtained
as $V\3 G$, where $V$ is flat and $G$ reductive.
Then $U$ is not isomorphic to $U'$ as
hyperk\"ahler manifold.

{\bf Proof:} $\dim_{tr}\Tw(U')=\dim \Tw(U')$ 
by \ref{_quotients_Moish_Theorem_}, and 
$\dim_{tr}\Tw(U)=1$ by \ref{_dim_tr_twistor_compact_Theorem_}.
\endproof

\hfill

A twistor space of a manifold obtained by hyperk\"ahler
reduction is Moishezon, and from the above argument
it is easy to see that it is Zariski open in
a compact Moishezon variety. However, it is (almost) never
quasiprojective. 

\hfill

\theorem\label{_twi_non_quasiproje_Theorem_}
Let $M$ be a hyperk\"ahler manifold such that its
twistor space $\Tw(M)$ can be embedded to a projective
manifold. Then, for each induced complex structure $L$, 
the complex manifold $(M,L)$ has no compact subvarieties
of positive dimension.

\hfill

{\bf Proof:} 
Consider the anticomplex involution $\iota$ on $\Tw(M)$ mapping
$(m,L)$ to $(m,-L)$.
Suppose that $(M,L)\subset \Tw(M)$ has a compact subvariety.
Since the $(M,L)$ is quasiprojective, this would imply
that $(M,L)$ contains a compact curve $S$. The curve $\iota(S)$
is also holomorphic in $(M,-L)$. 

Consider a K\"ahler form $\omega$ on $\Tw(M)$.
Since $\iota$ is antiholomorphic, 
$-\iota(\omega)$ is a closed, positive $(1,1)$-form.
Replacing $\omega$ by $\omega-\iota(\omega)$, 
we may assume that $\omega$ satisfies $\iota(\omega)=-\omega$.

Now, since the cohomology classes of $S$ and $\iota(S)$
are equal, we have 
\[ 
  \int_S \omega=\int_{\iota(S)}\omega=-\int_{\iota(S)}\iota(\omega)
\]
On the other hand, $\int_S \omega=\int_{\iota(S)}\iota(\omega)$
by functorial properties of integral. This implies
$\int_S \omega=0$, giving a contradiction.
\endproof

\hfill

\remark
As shown in \cite{_Verbitsky:subvar_},
a general fiber of the map $\Tw(M)\arrow \C P^1$
has only even-dimensional complex subvarieties.


\section{Appendix: Formal geometry and holography principle  (by Dmitry Kaledin)}
\label{_section_Appendix_}


\def\E{\mathcal{E}}
\def\F{\mathcal{F}}
\def\D{\mathcal{D}}
\def\C{\mathbb{C}}
\def\calo{\mathcal{O}}
\def\wh{\widehat}
\def\Hom{\operatorname{Hom}}
\def\RHom{\operatorname{RHom}}
\def\Hhom{\operatorname{\mathcal{H}om}}
\setlength{\unitlength}{1pt}
\newcommand{\hdot}{{\:\raisebox{3pt}{\text{\circle*{1.5}}}}}
\newcommand{\idot}{{\:\raisebox{1pt}{\text{\circle*{1.5}}}}}
\newcommand{\proof}[1][Proof.]{\smallskip\noindent{\em #1}}
\def\endproof{\hfill\ensuremath{\square}\par\medskip}

In this Appendix, we will try to explain the constructions of the
paper in a slightly more general context of $\D$-modules and
``formal geometry'' of Gelfand and Kazhdan \cite{_Gelfand_Kazhdan:formal_}. 
To save space,
we only sketch the proofs, and we work in the algebraic setting
(generalization to complex-analytic varieties is immediate, exactly
the same arguments work).

Assume given a smooth algebraic variety $X$ over a field $k$ of
characteristic $0$. We will work with coherent $\D$-modules over
$X$, that is, with sheaves of left modules over the algebra $\D_X$
of differential operators on $X$ which are finitely generated over
$\D$. Any coherent $\D$-module $\E$ is also a quasicoherent sheaf of
$\calo_X$-modules; recall that if $\E$ is coherent over $\calo_X$,
then it comes from a vector bundle on $X$ equipped with a flat
connection. In particular, the structure sheaf $\calo_X$ is a
$\D$-module; it corresponds to the trivial line bundle. For any
coherent sheaf $\E$ on $X$, we can consider the {\em induced
  $\D$-module} $\D_X \otimes_{\calo_X} \E$; this is coherent over
$\D_X$ but not over $\calo_X$. Coherent $\D$-modules form an abelian
category, and we can consider its derived category.

Assume given an open subvariety $U \subset X$ with the embedding map
$j:U \to X$, and let $Z = X \setminus U \subset X$ be the complement
$U$. Then the complex of quascoherent sheaves $R^\hdot j_*\calo_U$
on $X$ has a natural structure of a complex of $\D$-modules, and we
have an exact triangle
\begin{equation}\label{delta-eq}
\begin{CD}
\delta_Z @>>> \calo_X @>>> R^\hdot j_*\calo_X @>>>
\end{CD}
\end{equation}
of complexes of $\D$-modules on $X$, where $\delta_Z$ is supported
at $Z$. We will need the following standard result.

\bigskip

\lemma\label{delta-le}
In the notation above, for any coherent sheaf $\E$ on $X$, we have a
natural identification
$$
\Hom_{\calo_X}(\delta_Z,\E) \cong \Gamma(X,\E \otimes
\wh{\calo}_{X,Z}),
$$
where $\wh{\calo}_{X,Z}$ is the formal completion of structure sheaf
$\calo_X$ at the closed subscheme $Z \subset X$.

\medskip

\proof[Sketch of a proof.] Since $\delta_Z$ is supported at $Z$, the
natural map
$$
\Hom_{\calo_X}(\delta_Z,\E) \to \Hom_{\calo_X}(\delta_Z,\E \otimes
\wh{\calo}_{X,Z})
$$
is an isomorphism. On the other hand, we obviously have
$$
\Hom_{\calo_X}(R^\hdot j_*\calo_U,E \otimes \wh{\calo}_{X,Z}) = 0,
$$
and the claim then immediately follows from the long exact sequence
associated to the exact triangle \eqref{delta-eq}.
\endproof

For any map $f:X \to Y$ of smooth algebraic varieties, the pullback
functor $f^*$ extends to a functor between $\D$-modules, and its
derived functor extends to a functor between the derived categories
of $\D$-modules; we will denote this last functor by $f^?$. We have
$$
f^? = f^*[\dim X - \dim Y],
$$
where $f^*$ is the standard pullback functor for $\D$-modules, and
$[-]$ stands for cohomological shift. If $f$
is a smooth map, or a closed embedding, or a composition of the two,
then $f^?$ has a left-adjoint functor $f_?$ given by
$$
f_? = f_![\dim Y - \dim X],
$$
where $f_!$ is the standard functor of direct image with compact
supports. Then \ref{delta-le} has the following corollary.

\bigskip

\lemma\label{main-le}
Assume given a map $f:Y \to X$ of smooth algebraic varieties, and
assume that $f$ factors as
$$
\begin{CD}
Y @>{\pi}>> Z @>{\iota}>> X,
\end{CD}
$$
where $Z$ is smooth, $\iota$ is a closed embedding, and $\pi$ is a
smooth map with contractible fibers. Then for any vector bundle
$\E$ on $X$, we have a natural identification
\begin{equation}\label{main-eq}
\Hom_{\D_Y}(f^?(\D_X \otimes_{\calo_X} \E^*),\calo_Y) \cong
\Gamma(X,\E \otimes \wh{\calo}_{X,Z}),
\end{equation}
where $\E^*$ is the dual vector bundle, and $\D_X \otimes_{\calo_X}
\E^*$ is the corresponding induced $\D$-module.

\medskip

\proof{} For any $\D$-module $\F$ and $\calo_X$-module $\E$, we have
$$
\Hom_{\D_X}(\D_X \otimes_{\calo_X} \E,\F) \cong
\Hom_{\calo_X}(\E,F),
$$
so that \ref{delta-le} provides an identification
$$
\Hom_{\D_X}(\D_X \otimes_{\calo_X} \E^* \otimes_{\calo_X}
\delta_Z,\calo_X) \cong \Gamma(X,\E \otimes \wh{\calo}_{X,Z}),
$$
and by adjunction, it suffices to construct an isomorphism
$$
f_?f^?(\D_X \otimes_{\calo_X} \E^*) \cong \D_X \otimes_{\calo_X}
\E^* \otimes_{\calo_X} \delta_Z.
$$
But under our assumptions on $f$, we have $f_?\calo_X \cong
\delta_Z$, and we are done by the projection formula.
\endproof

We note that the left-hand side of \eqref{main-eq} admits a slightly
different interpretation. Recall that for any coherent sheaf $\E$ on
$X$, the {\em jet bundle} $J^\infty\E$ of $\E$ is a (pro)coherent
sheaf on $X$ given by
$$
J^\infty\E = \pi_{2*}\pi_1^*\E,
$$
where $\pi_1,\pi_2:\wh{X} \to X$ are the two natural projections of
the completion $\wh{X}$ of the product $X \times X$ near the
diagonal $X \subset X \times X$. In terms of $\D$-modules, we have
$$
J^\infty\E \cong \Hhom_{\calo_X}(\D_X,\E).
$$
The jet bundle $J^\infty\E$ carries a canonical flat connection, and
in the assumptions of \ref{main-le}, we have
$$
\begin{aligned}
\Hom_{\D_Y}(f^?(\D_X \otimes_{\calo_X} \E^*),\calo_Y) &\cong
\Hom_{\D_Y}(\calo_Y,\Hhom_{\calo_Y}(f^?(\D_X \otimes_{\calo_X}
\E^*),\calo_Y))\\
&\cong \Hom_{\D_Y}(\calo_Y,f^?(\Hhom_{\calo_X}(\D_X \otimes_{\calo_X}
\E^*,\calo_X)))\\
& \cong \Gamma^\nabla(Y,f^*J^\infty\E),
\end{aligned}
$$
where $\Gamma^\nabla(-)$ stands for the space of flat global
sections. Then \eqref{main-eq} reads as
\begin{equation}\label{jet-eq}
\Gamma^\nabla(Y,f^*J^\infty\E) \cong \Gamma(X,\E \otimes
\wh{\calo}_{X,Z}).
\end{equation}

\bigskip

\remark
The assumptions of \ref{main-le} are in fact too strong. Firstly, it
is clearly enough to require that the fibers of $\pi$ are non-empty
and connected, so that their top degree cohomology with compact
supports is one-dimensional -- and under the assumptions as stated,
we not only obtain an isomorphism of $\Hom$'s but also of the
$\RHom$'s, so that \eqref{jet-eq} extends to an isomorphism
$$
H^\hdot_{DR}(Y,f^*J^\infty\E) \cong H^\hdot(X,\E \otimes
\wh{\calo}_{X,Z})
$$
of cohomology groups. Secondly, one probably does not have to
require that $\pi$ and $Z$ are smooth -- some assumptions are
needed, but they can considerably relaxed. However, since even the
stronger assumptions work for us, we did not pursue this.

\bigskip

Assume now that we are given smooth proper algebraic varieties $Y$,
$X$, and a family of closed embedding from $Y$ to $X$ parametrized
by a smooth algebraic variety $T$ -- that is, we have a map
$$
f:Y \times T \to X
$$
such that for any $t \in T$, the corresponding map $f_t:Y = Y \times
t \to X$ is a closed embedding. Moreover, assume given a coherent
sheaf $\E$ on $X$, denote by $\rho:Y \times T \to T$ the projection,
and denote
$$
\Phi_f\E = \rho_*^\nabla f^*J^\infty\E,
$$
where $f^*J^\infty\E$ is the pullback of the jet bundle $J^\infty\E$
equipped with its natural flat connection $\nabla$, and
$\rho_*^\nabla$ stands for sheaf of relative flat sections. Then
$\Phi_a\E$ is a sheaf on $T$, the base of the family, and by
\eqref{jet-eq}, the fiber $(\Phi_f\E)_t$ at a point $t \in T$ is
given by
$$
(\Phi_f\E)_t = \Gamma^\nabla(Y,f_t^*J^\infty\E) \cong \Gamma(X,\E
\otimes \wh{\calo}_{X,Y_t}),
$$
where $Y_t \subset X$ is the image of the closed embedding $f_t:Y
\to X$. On the other hand, $\Phi_f\E$ carries a natural flat
connection, and if we assume that $f$ is smooth with contractible
fibers, we have
$$
\Gamma^\nabla(T,\Phi_f\E) = \Gamma^\nabla(T \times Y,f^*J^\infty\E)
\cong \Gamma(X,\E),
$$
again by \eqref{jet-eq}.

We now note that this is exactly the situation that we have in the
paper. Namely, we take $X$ to be the twistor space of a
hyperk\"ahler manifold $M$, we take $T = M_{\C}$ to be the
complexification of the real-analytic manifold underlying $M$, we
take $Y = \C P^1$, and we let
$$
f:T \times Y \to X
$$
be the standard family of twistor lines (real points in this family
correspond to horizontal sections of the twistor fibration $X \to \C
P^1$ parametrized by points of $M$). Then since the normal bundle to
any line in our family is a sum of several copies of $\calo(1)$, the
family is unobstructed, and moreover, it remains unobstructed even
if we fix a point at a twistor line, so that the map $f$ is
smooth. Its fibers are small polydiscs, thus contractible, and all
the assumptions of \ref{main-le} are therefore satisfied.

\hfill

{\small

\hfill

\noindent {\sc Misha Verbitsky\\
Laboratory of Algebraic Geometry, \\
Faculty of Mathematics, National Research University HSE,\\
7 Vavilova Str. Moscow, Russia
 }

}


\begin{thebibliography}{MM}

\bibitem[BBI]{_Badescu_Beltrametti_Ionescu_}
B\v adescu, Lucian; Beltrametti, Mauro C.; Ionescu, Paltin,
{\em Almost-lines and quasi-lines on projective manifolds}, 
Complex analysis and algebraic geometry, 1-27, de Gruyter, Berlin, 2000. 

\bibitem[AHS]{_Atiyah_Hitchin_Singer_}
Atiyah, M. F., Hitchin, N. J., and Singer, I. M. {\em Self-duality in
four-dimensional Riemannian geometry}, Proc. Roy. Soc. London Ser. A
362 (1978), 425-461.


\bibitem[Bes]{_Besse:Einst_Manifo_}
Besse, A., {\em Einstein Manifolds}, Springer-Verlag, New York (1987).

\bibitem[C]{_Campana:twistor_}
F. Campana, {\em On twistor spaces of the class ${\cal C}$}, 
J. Differential Geom. 33 (1991) 541-549.


\bibitem[D]{_Demailly:pseudoconvex-concave_}
Demailly, Jean-Pierre,
{\em Pseudoconvex-concave duality and regularization of currents,} 
Several complex variables (Berkeley, CA, 1995-1996), 233-271,
Math. Sci. Res. Inst. Publ., 37, Cambridge Univ. Press, Cambridge, 1999. 

\bibitem[GK]{_Gelfand_Kazhdan:formal_}
I.M. Gelfand and D.A. Kazhdan, {\em Some problems
of differential geometry and the calculation of cohomologies of Lie
algebras of vector fields}, Soviet Math. Dokl. {\bf 12} (1971),
1367-1370.


\bibitem[Har]{_Hartshorne:cohomo-1968_}
Hartshorne, Robin,
{\em Cohomological dimension of algebraic varieties,}
Ann. of Math. (2) 88 1968 403-450. 

\bibitem[HKLR]{_HKLR_}  
N. J. Hitchin, A. Karlhede, U. Lindstr\"om, M. Ro\v cek, 
{\em Hyperk\"ahler metrics and supersymmetry,}
Comm. Math. Phys. {\bf 108}, (1987) 535--589.

\bibitem[Ho]{_Honda:classification_}
Nobuhiro Honda
{\em Moishezon twistor spaces on $4\C P^2$},
arXiv:1112.3109, 51 pages.

\bibitem[K]{_Kaledin:twistor_}
D. Kaledin,
{\em Integrability of the twistor space for a hypercomplex manifold},
Sel. math., New ser. {\bf 4} (1998) 271-278.


\bibitem[KV]{_NHYM_} 
Kaledin, D., Verbitsky, M.,
{\it Non-Hermitian 
Yang-Mills connections}, Selecta Math. (N.S.) {\bf 4}
(1998), no. 2, 279--320.

\bibitem[KST]{_Kebekus_Sola_Toma_}
Stefan Kebekus, Luis Sola Conde, Matei Toma,
{\em 
Rationally connected foliations after Bogomolov and McQuillan,}
J. Algebraic Geom. 16 (2007), no. 1, 65-81. 

\bibitem[Ko]{_Kollar:curves_}
Koll\'ar, J., 
{\em Rational curves on algebraic varieties,} Springer, 
  1996, viii+320 pp..


\bibitem[Nak]{_Nakajima_}
H. Nakajima,
{\em Lectures on Hilbert schemes of points on surfaces},
Providence: American Mathematical Society, 1999.

\bibitem[P]{_Poon:twistors_}
Y. S. Poon, {\em 
On the algebraic structure of twistor spaces,} 
J. Diff. Geom. 36 (1992), 451-491.

\bibitem[SV]{_Soldatenkov_Verbitsky:triana_}
Andrey Soldatenkov, Misha Verbitsky,
{\em Subvarieties of hypercomplex manifolds with holonomy in $SL(n,{\Bbb H})$},
Journal of Geometry and Physics, Volume 62, Issue 11
(2012), Pages 2234-2240, arXiv:1202.0222

\bibitem[V1]{_Verbitsky:trianalyt_}
Verbitsky M., {\em Tri-analytic subvarieties of hyper-Kaehler manifolds,} 
also known as {\em Hyperk\"ahler embeddings and holomorphic 
symplectic geometry II}, GAFA {\bf 5} no. 1 (1995), 92-104, 
alg-geom/9403006.


\bibitem[V2]{_Verbitsky:hypercomple_}
Verbitsky, M.,
{\em Hypercomplex Varieties}, 
alg-geom/9703016,
 Comm. Anal. Geom. {\bf 7} 
(1999), no. 2, 355--396.


\bibitem[V3]{_Verbitsky:subvar_}
Verbitsky, M., {\em Subvarieties in non-compact hyperk\"ahler manifolds},
Math. Res. Lett., vol. 11 (2004), no. 4, pp. 413--418.

\bibitem[V4]{_Verbitsky:rational_twi_}
Verbitsky, M., {\em Rational curves and special metrics on
  twistor spaces}, arXiv:1210.6725, 12 pages.


\bibitem[Y]{_Yau:Calabi-Yau_} 
Yau, S. T., {\em On the Ricci curvature of a compact K\"ahler manifold 
and the complex Monge-Amp\`ere equation I.}  Comm. on Pure and Appl.
Math. 31, 339-411 (1978).


\end{thebibliography}
\end{document}